\documentclass{amsart}

\usepackage{amsmath,amssymb,amsthm}
\usepackage{graphicx,tikz, subcaption}

\newtheorem*{thm}{Theorem}

\theoremstyle{definition}

\theoremstyle{remark}

\begin{document}

\title[]{Free Convolution Powers \\via Roots of Polynomials}
\thanks{S.S. is supported by the NSF (DMS-1763179) and the Alfred P. Sloan Foundation.}

\author[]{Stefan Steinerberger}
\address{Department of Mathematics, University of Washington, Seattle}
\email{steinerb@uw.edu}

\begin{abstract} Let $\mu$ be a compactly supported probability measure on the real line. Bercovici-Voiculescu and Nica-Speicher proved the existence of a free convolution power $\mu^{\boxplus k}$ for any real $k \geq 1$. The purpose of this short note is to give a formal proof of an elementary description of $\mu^{\boxplus k}$ in terms of of polynomials and roots of their derivatives. This bridge allows us to switch back and forth between free probability and the asymptotic behavior of polynomials. \end{abstract}

\maketitle

\section{Introduction}
\subsection{Free Convolution.} The notion of free convolution $\mu \boxplus \nu$ of two compactly supported probability measures is due to Voiculescu \cite{voic}. A definition is as follows: for any compactly supported probability measure, we can consider its Cauchy transform $G_{\mu}:\mathbb{C} \setminus \mbox{supp}(\mu) \rightarrow \mathbb{C}$ defined via
$$ G_{\mu}(z) = \int_{\mathbb{R}} \frac{1}{z-x} d\mu(x).$$
For a compactly supported measure, $G_{\mu}(z)$ tends to 0 as $|z| \rightarrow \infty$.
Given $G_{\mu}$, we define the $R-$transform $R_{\mu}(s)$ for sufficiently small complex $s$ by demanding that
$$ \frac{1}{G_{\mu}(z)} + R_{\mu}(G_{\mu}(z)) = z$$
for all sufficiently large $z$. The free convolution $\mu \boxplus \nu$ is then the unique compactly supported measure for which
$$ R_{\mu \boxplus \nu}(s) = R_{\mu}(s) + R_{\nu}(s)$$
for all sufficiently small $s$. A fundamental result due to Voiculescu is the \textit{free central limit theorem}: if $\mu$ is a compactly supported probability measure with mean 0 and variance 1, then suitably rescaled copies of $\mu^{\boxplus k}$ converge to the semicircular distribution. This notion can be extended to real powers.

\begin{thm}[Fractional Free Convolution Powers exist, \cite{berc, nica}] Let $\mu$ be a compactly supported probability measure on $\mathbb{R}$ and assume $k \geq 1$ is real. Then there exists a unique compactly supported probability measure $\mu^{\boxplus k}$ such that
$$ R_{\mu^{\boxplus k}}(s) = k \cdot R_{\mu}(s) \qquad \quad \mbox{for all $s$ sufficiently small.}$$
\end{thm}
This was first shown for $k$ sufficiently large by Bercovici \& Voiculescu \cite{berc} and then by Nica \& Speicher \cite{nica} for all $k \geq 1$. We also refer to \cite{ans, bel0, bel1, bel2, berc1, hiai, huang, mingo, nica2, shlya, williams}. The purpose of this note is to (formally) prove an elementary description of $\mu^{\boxplus k}$ in terms of polynomials and the density of the roots of their derivatives.

\subsection{Polynomials.} Roots of polynomials are a classical subject and there are many results we do not describe here, see \cite{byun, farmer, gauss, granero, han, huang, kab, kab2, kornik, lucas, nica, or, or3, or4, or2, pem, pem2, pol, ravi, red, sub, steini, riesz}. Our problem will be as follows:  let $\mu$ be a compactly supported probability measure on the real line and suppose $x_1, \dots, x_n$ are $n$ independent random variables sampled from $\mu$ (which we assume to be sufficiently nice).
We then associate to these numbers the random polynomial
$$ p_n(x) = \prod_{k=1}^{n}{(x-x_k)}$$
having roots exactly in these points. What can we say about the behavior of the roots of the derivative $p_n'$? There is an interlacing phenomenon and the roots of $p_n$ are also distributed according to $\mu$ as $n \rightarrow \infty$. The same is true for the second derivative $p_n''$ and any finite derivative. However, once the number of derivatives is proportional to the degree, the distribution will necessarily change.
\begin{quote}
\textbf{Question.} Fix $0 < t < 1$. How are the roots of $p_n^{\left\lfloor t \cdot n\right\rfloor}$ distributed?
\end{quote}
The question was raised by the author in \cite{steini}. 
The answer, if it exists, should be another measure $u(t,x)dx$. Note that, since this measure describes the distribution of roots of polynomials of degree $(1-t) \cdot n$, as $n \rightarrow \infty$, we have
$$ \int_{\mathbb{R}} u(t,x) dx = 1-t.$$
Relatively little is known about the evolution of $u(t,x)$: \cite{steini} established, on a formal level, a PDE for $u(t,x)$. This PDE is given by
$$ \frac{\partial u }{\partial t} + \frac{1}{\pi} \frac{\partial}{\partial x} \arctan\left( \frac{Hu}{u} \right) = 0 \qquad \mbox{on}~\mbox{supp}(u),$$
where
$$ Hf(x) =  \mbox{p.v.}\frac{1}{\pi} \int_{\mathbb{R}}{\frac{f(y)}{x-y} dy} \qquad \mbox{is the Hilbert transform.}$$

\begin{figure}[h!]
\begin{minipage}[l]{.45\textwidth}
\begin{tikzpicture}
\node at (0,0) {\includegraphics[width = 0.7\textwidth]{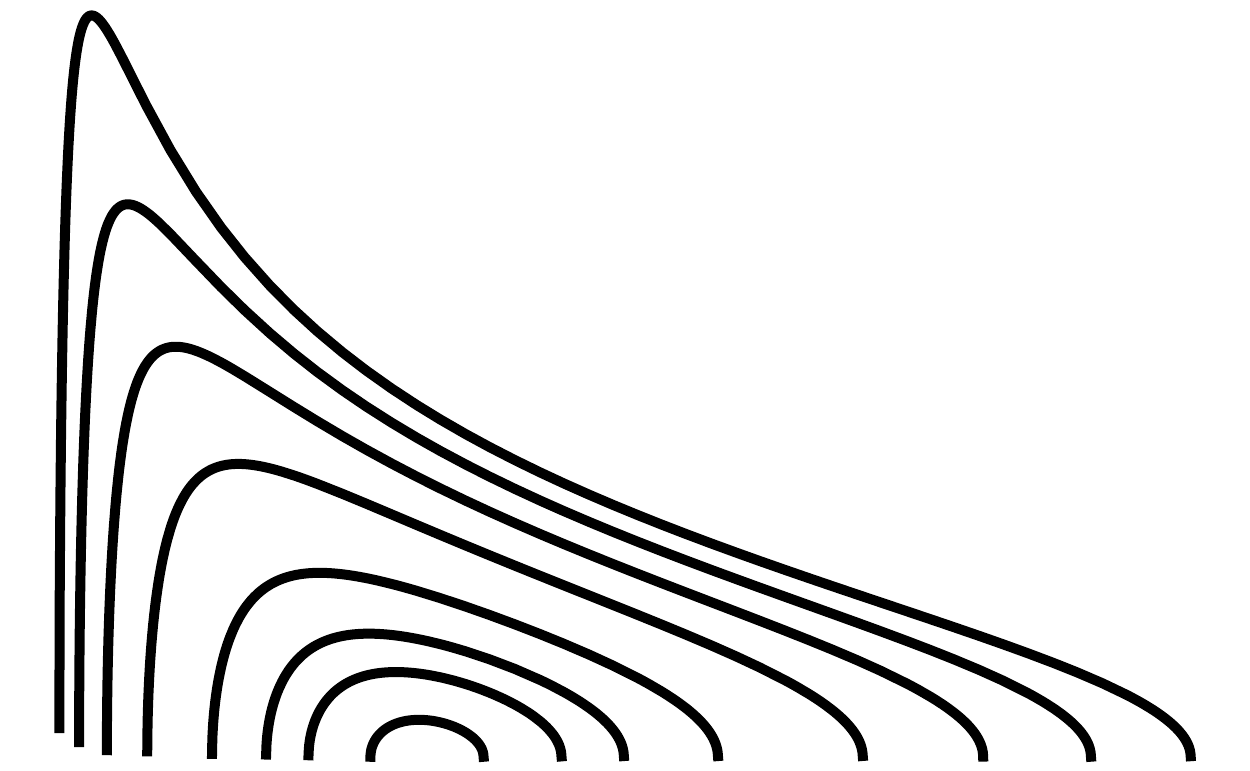}};
\end{tikzpicture}
\end{minipage} 
\begin{minipage}[r]{.45\textwidth}
\begin{tikzpicture}
\node at (0,0) {\includegraphics[width = 0.7\textwidth]{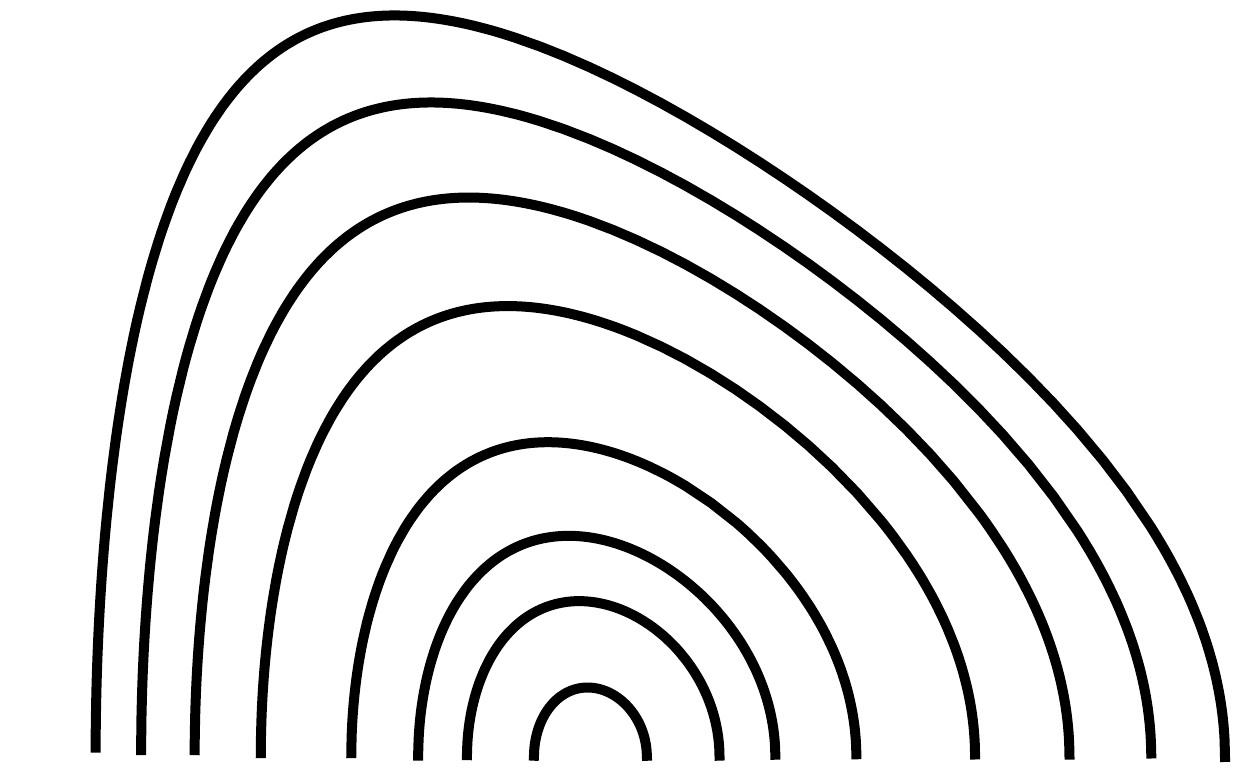}};
\end{tikzpicture}
\end{minipage} 
\caption{The densities of two evolving measures $u(t,x)dx$. They shrink and vanish at time $t=1$.} 
\end{figure}

It is also known that it has to satisfy the conservation laws
$$\int_{\mathbb{R}} \int_{\mathbb{R}} u(t,x) (x-y)^2 u(t,y) ~ dx dy = (1-t)^3 \int_{\mathbb{R}} \int_{\mathbb{R}} u(0,x) (x-y)^2 u(0,y) ~ dx dy.$$
 Hoskins and the authors \cite{hoskins} established a universality result for large derivatives of polynomials with random roots: such derivatives behave like random shifts of Hermite polynomials. Hermite polynomials, in turn, have roots whose density is given by a semicircle and this leads one to believe that $u(t,x)$ should, for $t$ close to 1, look roughly like a semicircle (and this has also been observed numerically).

There are two explicit closed-form solutions, derived in \cite{steini}, a shrinking semicircle and a one-parameter solution that lies in the Marchenko-Pastur family (see Fig. 1).
Numerical simulations in \cite{hoskins} also suggested that the solution tends to become smoother. O'Rourke and the author \cite{or2} derived an analogous transport equation
for polynomials with roots following a radial distribution in the complex plane.

\section{The Result}
\subsection{An Equivalence.}
We can now state our main observation: both the free convolution of a measure with itself, $\mu^{\boxplus k}$, and the density of roots of derivatives of polynomials, $u(t,x)$, are described by the same underlying process. 
\begin{thm} At least formally, if $\mu = u(0,x)dx$ and $\left\{x : u(0,x) > 0\right\}$ is an interval, then for all real $k\geq 1$
$$ \mu^{\boxplus k} =u\left(1- \frac{1}{k}, \frac{x}{k} \right) dx.$$
\end{thm}
We first clarify the meaning of `formally'. In a recent paper, Shlyakhtenko \& Tao \cite{shlya} derived, formally, a PDE for the evolution of the $\mu^{\boxplus k}$. This PDE happens to be the same PDE (expressed in a different coordinate system) that was formally derived by the author for the evolution of $u(t,x)$ \cite{steini}. The derivation in \cite{steini} is via a mean-field limit approach, the `microscopic' derivation is missing. In particular, the derivation in \cite{steini} assumed the existence of $u(t,x)$ and a crystallization phenomenon for the roots; such a crystallization phenomenon has been conjectured for a while, there is recent progress by Gorin \& Kleptsyn \cite{gorin}.\\

    \begin{figure}[h!]
        \includegraphics[width=0.6\textwidth]{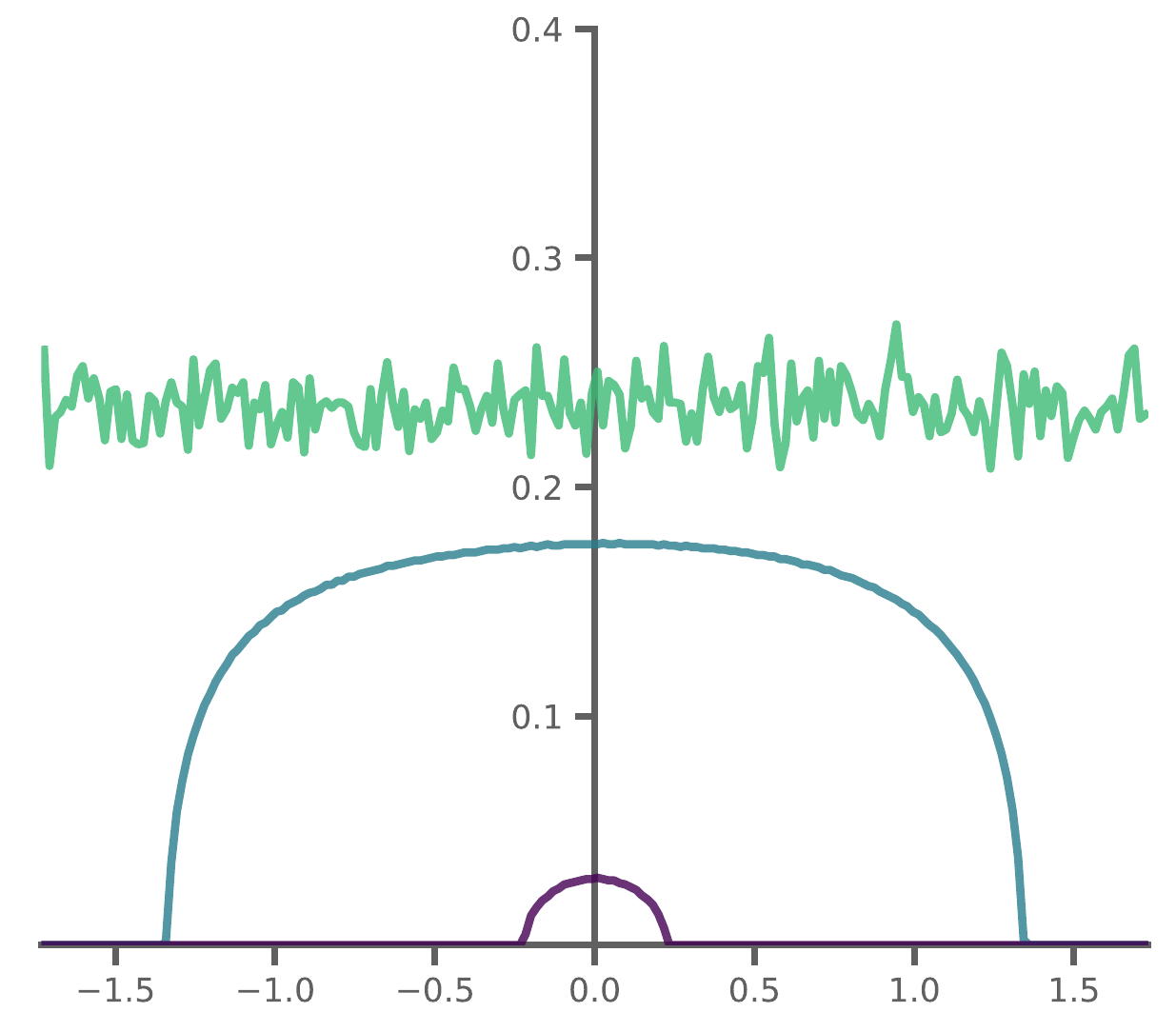}
\caption{Evolution of $u(t,x)$ (from \cite{hoskins}) starting with random and uniformly distributed roots: the evolution smoothes and we see a semicircle before it vanishes.}
\end{figure}

Naturally, this has a large number of consequences since it allows us to go back and forth between results from
free probability and results regarding polynomials and their roots. 
As an illustration, we recall that for the semicircle law $\mu_{sc}^{\boxplus k}$ is a another semicircle law stretched by a factor $k^{1/2}$, thus
$$ \mu_{sc}^{\boxplus k} = \frac{2}{\pi} \sqrt{ \frac{1}{k} - \frac{x^2}{k^2}}dx$$
Conversely, as was computed in \cite{steini}, the evolution of densities of polynomials when beginning with a semicircle behaves as
$$ u(t,x) = \frac{2}{\pi}\sqrt{1-t- x^2}.$$

One immediate consequence of the equivalence is that it provides us with a fast algorithm to approximate $\mu^{\boxplus k}$ when $\mu = f(x)dx$ and $f$ is smooth. This may be useful in the study of the semigroup $\mu^{\boxplus k}$. Using the logarithmic derivative $p_n'/p_n$, it is possible quickly differentiate real-rooted polynomials $p_n$ a large number of times, $t \cdot n$, even when the degree is as big as $n \sim 100.000$: this was done in \cite{hoskins} using a multipole method (a modification of an algorithm due to Gimbutas, Marshall \& Rokhlin \cite{gimb}). Fig. 2 shows an example computed using 80.000 roots: we observe the initial smoothing and the eventual convergence to a semicircle.

\subsection{Some Connections.} Some connections are as follows.\\

\textit{The Free Central Limit Theorem.} Voiculescu \cite{voic} proved that $\mu^{\boxplus k}$ (suitably rescaled) approaches a semicircle distribution in the limit.  Motivated by high-precision numerics, Hoskins and the author \cite{hoskins} conjectured that $u(t,x)$ starts looking like a semicircle for $t$ close to 1 and proved a corresponding universality result for polynomials with random roots: if $p_n$ is a polynomial with random roots (from a probability measure $\mu$ whose moments are all finite), then, for fixed $\ell \in \mathbb{N}$ and $n \rightarrow \infty$, we have for $x$ in a compact interval
 $${n^{\ell/2}} \frac{\ell!}{n!} \cdot p_n^{(n-\ell)}\left( \frac{x}{\sqrt{n}}\right) \rightarrow He_{\ell}(x + \gamma_n),$$
 where $He_{\ell}$ is the $\ell-$th probabilists' Hermite polynomial and $\gamma_n$ is a random variable converging to the standard $\mathcal{N}(0,1)$ Gaussian as $n \rightarrow \infty$. Hermite polynomials have roots that are asymptotically distributed like a semicircle. A result in the deterministic setting has recently been provided by Gorin \& Kleptsyn \cite{gorin}. \\
 
  \textit{Conservation Laws.} The author showed that the evolution $u(t,x)$ satisfies the algebraic relations 
\begin{align*}
 \int_{\mathbb{R}}{ u(t,x) ~ dx} = 1-t, \qquad \qquad  \int_{\mathbb{R}}{ u(t,x) x ~ dx} = \left(1-t\right)\int_{\mathbb{R}}{ u(0,x) x~ dx}, \qquad\\
  \int_{\mathbb{R}} \int_{\mathbb{R}} u(t,x) (x-y)^2 u(t,y) ~ dx dy = (1-t)^3 \int_{\mathbb{R}} \int_{\mathbb{R}} u(0,x) (x-y)^2 u(0,y) ~ dx dy.
\end{align*}
These are derived from Vieta-type formulas that express elementary symmetric polynomials in terms of power sums. Equivalently, we have
$ \kappa_n ( \mu^{\boxplus k}) = k^n \kappa_n(\mu),$
where $\kappa_n$ is the $n-$th free cumulant providing a large number of conservation laws.\\

\textit{Monotone Quantities.}
Voiculescu \cite{voic2} introduced the free entropy
$$ \chi(\mu) = \int_{\mathbb{R}} \int_{\mathbb{R}} \log{ |s-t|} d\mu(s) d\mu(t) + \frac{3}{4} + \frac{\log{(2\pi)}}{2}$$
and the free Fisher information
$$\Phi(\mu) = \frac{2\pi^2}{3} \int_{\mathbb{R}} \left( \frac{d\mu}{dx} \right)^3 dx.$$ 
Shlyakhtenko \cite{shlya0} proved that $\chi$ increases along free convolution of $\mu$ with itself whereas $\Phi$ decreases (both suitably rescaled). Shlyakhtenko \& Tao \cite{shlya} showed monotonicity along the entire flow $\mu^{\boxplus k}$ for real $k \geq 1$. Conversely, on the side of polynomials, it is known that
$$ \frac{\left|\left\{x \in \mathbb{R}: u(t,x) > 0 \right\}\right|}{1-t} \qquad \mbox{is non-decreasing in time.}$$
 Another basic result for polynomials is commonly attributed to Riesz \cite{farmer, riesz}: denoting the smallest gap of a polynomial $p_n$ having $n$ real roots $\left\{x_1, \dots, x_n\right\}$ by 
$$ G(p_n) = \min_{i \neq j}{|x_i - x_j|},$$
we have $G(p_n') \geq G(p_n)$: the minimum gap grows under differentiation. A simple proof is given by Farmer \& Rhoades \cite{farmer}. This would suggest that
the maximal density cannot increase over time.\\

\textit{The Minor Process.} Shlyakhtenko \& Tao \cite{shlya} connect the evolution to the minor process: trying to understand how the eigenvalues of the $n \times n$ minor of a large random Hermitian matrix $N \times N$ behave. This answers a question numerically verified by Hoskins and the author \cite{hoskins}.\\

\textit{Related Results.} There are several other papers in the literature that seem to be connected to this circle of ideas. We mention Gorin \& Marcus \cite{gorin0}, Marcus \cite{marcus0}, Marcus, Spielman \& Srivastava \cite{marcus}.

\section{Proof}
\begin{proof} Shlyakhtenko \& Tao \cite{shlya} derive that if
$$ d\mu^{\boxplus k} = f_k(x) dx$$
and if we substitute $k=1/s$ (thus $0< s < 1$) and $f := f_{1/s}$, then on a purely formal level
$$ \left(-s \frac{\partial}{\partial s} + x \frac{\partial}{ \partial x}\right) f = \frac{1}{\pi} \frac{\partial}{\partial x} \arctan \left( \frac{f}{Hf} \right).$$
On the other hand, the author derived \cite{steini}, also on a formal level, that as long as $\left\{x: u(t,x) > 0\right\}$ is an interval
$$ \frac{\partial u }{\partial t} + \frac{1}{\pi} \frac{\partial}{\partial x} \arctan\left( \frac{Hu}{u} \right) = 0 \qquad \mbox{on}~\mbox{supp}(u).$$
We note that Huang \cite{huang} showed that the number of connected components in the support of $\mu^{\boxplus k}$ is non-decreasing in $k$ which shows that once the support is an interval, this property is preserved.
We want to show that the solutions of these two PDEs are related via a change of variables: since both evolutions obey the same PDE, they must coincide. We observe that one nonlinear term seems to be
the reciprocal of the other, however, this compensates for the different sign.
We compute
\begin{align*}
\frac{\partial}{\partial x} \arctan\left( \frac{f}{Hf} \right) = \frac{1}{1+\frac{f^2}{(Hf)^2}} \partial_x \frac{f}{Hf} =\frac{ f_x (Hf) - f (Hf)_x }{f^2 + (Hf)^2}
\end{align*}
and compare it to
\begin{align*}
\frac{\partial}{\partial x}  \arctan\left( \frac{Hf}{f} \right) = \frac{1}{1+\frac{(Hf)^2}{f^2}} \partial_x \frac{Hf}{f} = \frac{  (Hf)_x f - f_x (Hf) }{f^2 + (Hf)^2}
\end{align*}
and see that it is the same term with opposite sign. This allows us to write
$$ \frac{\partial u }{\partial t} =  \frac{1}{\pi} \frac{\partial}{\partial x} \arctan\left( \frac{u}{Hu} \right).$$
We now claim that
$$ f(s,x) = u(1-s, s x).$$
Note that the left-hand side transforms
\begin{align*}
 (-s \partial_s + x \partial_x) f = -s \left( \frac{\partial u}{\partial t} (-1) + \frac{\partial u}{\partial x} x\right) + x \frac{\partial u}{\partial x} s = s \frac{\partial u}{\partial t}.
\end{align*}
It remains to understand how the right-hand side transforms. The Hilbert transform commutes with dilations and thus
$$\arctan\left( \frac{f}{Hf} \right)= \arctan\left( \frac{u(1-s, sx)}{H \left[ u(1-s, sx)\right]} \right)= \arctan\left(  \frac{u(1-s, sx)}{\left[H u (1-s, \cdot) \right] (sx)}\right)$$ 
whose derivative scales exactly by a factor of $s$.
\end{proof} 

\textbf{Remarks.} We see that, both derivations being purely formal, many problems remain. Indeed, this connection suggests
many interesting further avenues to pursue. 
 Roots of polynomials seem to regularize under differentiation at the micro-scale: if one
were to take a polynomial with random (or just relatively evenly spaced roots), then the roots of the $(\varepsilon \cdot n)-$th derivative are conjectured to behave
locally like arithmetic progressions up to a small error. Results of this flavor date back to Polya \cite{pol} for analytic functions, see also Farmer \& Rhoades \cite{farmer}
and Pemantle \& Subramanian \cite{pem2}.  In the converse direction, it could be interesting to study the behavior of $u(t,x)$ when the initial conditions are
close to a semi-circle: despite the equation being both non-linear and non-local, its linearization around the semicircle seems to diagonalize nicely under
Chebychev polynomials -- can PDE techniques be used to get convergence rates for the free central limit theorem?

\end{document}